\documentclass{amsart}
\usepackage[dvips]{graphicx}
\usepackage{amscd}
\usepackage{amsmath}
\usepackage{amsxtra}
\usepackage{amsfonts}
\usepackage{amssymb}

\newtheorem{theorem}{Theorem}[section]
\newtheorem{corollary}[theorem]{Corollary}
\newtheorem{lemma}[theorem]{Lemma}
\newtheorem{proposition}[theorem]{Proposition}
\theoremstyle{definition}
\newtheorem{definition}[theorem]{Definition}
\newtheorem{remark}[theorem]{Remark}

\theoremstyle{remark}

\renewcommand{\theclaim}{\textup{\theclaim}}

\renewcommand{\theenumi}{\roman{enumi}}

\numberwithin{equation}{section}

\def\openone

{\mathchoice

{\hbox{\upshape \small1\kern-3.3pt\normalsize1}}

{\hbox{\upshape \small1\kern-3.3pt\normalsize1}}

{\hbox{\upshape \tiny1\kern-2.3pt\SMALL1}}

{\hbox{\upshape \Tiny1\kern-2pt\tiny1}}}

\makeatletter

\newbox\ipbox

\newcommand{\ip}[2]{\left\langle #1\,|\,#2\right\rangle}

\newcommand{\diracb}[1]{\left\langle #1\mathrel{\mathchoice

{\setbox\ipbox=\hbox{$\displaystyle \left\langle\mathstrut
#1\right.$}

\vrule height\ht\ipbox width0.25pt depth\dp\ipbox}

{\setbox\ipbox=\hbox{$\textstyle \left\langle\mathstrut
#1\right.$}

\vrule height\ht\ipbox width0.25pt depth\dp\ipbox}

{\setbox\ipbox=\hbox{$\scriptstyle \left\langle\mathstrut
#1\right.$}

\vrule height\ht\ipbox width0.25pt depth\dp\ipbox}

{\setbox\ipbox=\hbox{$\scriptscriptstyle \left\langle\mathstrut
#1\right.$}

\vrule height\ht\ipbox width0.25pt depth\dp\ipbox}

}\right. }

\newcommand{\dirack}[1]{\left. \mathrel{\mathchoice

{\setbox\ipbox=\hbox{$\displaystyle \left.\mathstrut
#1\right\rangle$}

\vrule height\ht\ipbox width0.25pt depth\dp\ipbox}

{\setbox\ipbox=\hbox{$\textstyle \left.\mathstrut
#1\right\rangle$}

\vrule height\ht\ipbox width0.25pt depth\dp\ipbox}

{\setbox\ipbox=\hbox{$\scriptstyle \left.\mathstrut
#1\right\rangle$}

\vrule height\ht\ipbox width0.25pt depth\dp\ipbox}

{\setbox\ipbox=\hbox{$\scriptscriptstyle \left.\mathstrut
#1\right\rangle$}

\vrule height\ht\ipbox width0.25pt depth\dp\ipbox}

} #1\right\rangle}

\newcommand{\lonet}{L^{1}\left(  \mathbb{T}\right)}

\newcommand{\linft}{L^{\infty}\left(  \mathbb{T}\right)}

\newcommand{\ltwor}{L^{2}\left(\mathbb{R}\right)}

\newcommand{\linfr}{L^{\infty}\left(\mathbb{R}\right)}

\input cyracc.def

\begin{document}
\title[The Wavelet Galerkin Operator ]{The Wavelet Galerkin Operator }
\author{Dorin E. Dutkay }
\email{ddutkay@math.uiowa.edu}
\thanks{}
\subjclass{}
\keywords{}
\begin{abstract}
We consider the eigenvalue problem
$$R_{m_0,m_0}h=\lambda h,\quad h\in C(\mathbb{T}),\quad |\lambda|=1$$
where $R_{m_0,m_0}$ is the wavelet Galerkin operator associated to
a wavelet filter $m_0$. The solution involves the construction of
representations of the algebra $\mathfrak{A}_N$- the $C^*$-algebra
generated by two unitaries $U,V$ satisfying $UVU^{-1}=V^N$
introduced in \cite{Jor}.
\end{abstract}
\maketitle
\tableofcontents
\section{\label{intro} Introduction}
 The wavelet Galerkin operator appears in several different contexts
such as wavelets (see for example \cite{Law91b}, \cite{Dau95},
\cite{CoDa96}, \cite{LLS96}, \cite{CoRa90}), ergodic theory and
$g$-measures (\cite{Kea72}) or quantum statistical mechanics
(\cite{Rue68}). For some of the applications of the Ruelle
operator we refer the reader to the book by V. Baladi
\cite{Bal00}. It also bears many different names in the
literature: the Ruelle operator, the Perron-Frobenius-Ruelle
operator, the Ruelle-Araki operator, the Sinai-Bowen-Ruelle
operator, the transfer operator and several others. We used the
name wavelet Galerkin operator as suggested in \cite{Law91b},
because of its close connection to wavelets that we will be using
in the sequel. We will also use the name Ruelle operator and
transfer operator.
\par
 The Ruelle operator considered in this paper is defined by
$$R_{m_0,m_0'}f(z)=\frac{1}{N}\sum_{w^N=z}\overline{m_0(w)}m_0'(w)f(w),\quad(z\in\mathbb{T})$$
where $m_0,m_0'\in\linft$ are nonsingular (i.e. they do not vanish
on a set of positive measure ), $\mathbb{T}$ is the unit circle
$\{z\in\mathbb{C}\,|\,|z|=1\}$, $N\geq 2$ is an integer. A large
amount of information about this operator is contained in
\cite{BraJo}. One of the main objectives of this paper is to do a
peripheral spectral analysis for the Ruelle operator, that is to
solve the equation
$$R_{m_0,m_0}h=\lambda h,\quad |\lambda|=1,\quad h\in C(\mathbb{T}).$$
\par
The restrictions that we will impose on $m_0$ are:
\begin{equation}\label{relm01}
m_0\in\mbox{Lip}_1(\mathbb{T})
\end{equation}
where
$$\mbox{Lip}_1(\mathbb{T})=\left\{ f:\mathbb{T}\rightarrow\mathbb{C}\,|\, f \mbox{ is Lipschitz }\right\}.$$
\begin{equation}\label{relm02}
m_0\mbox{ has a finite number of zeros }.
\end{equation}
\begin{equation}\label{relm03}
R_{m_0,m_0}1=1
\end{equation}
\begin{equation}\label{relm04}
m_0(1)=\sqrt{N}
\end{equation}

\par
In ergodic theory the Ruelle operators are used in the derivation
of correlation inequalities (see \cite{Sto01} and \cite{Dol98})
and in understanding the Gibbs measures. The role played by the
Ruelle operator in wavelet theory is somewhat similar. It can be
used to make a direct connection to the cascade approximation and
orthogonality relations.
\par
In the applications to wavelets, the function $m_0$ is a wavelet
filter, i.e., its Fourier expansion
\begin{equation}\label{eqj1}
m_0(z)=\sum_{k\in\mathbb{Z}}a_kz^k
\end{equation}
yields the masking coefficients of the scaling function $\varphi$
on $\mathbb{R}$, i.e. the function which results from the the
fixed-point problem
\begin{equation}\label{eqj2}
\varphi(x)=\sqrt{N}\sum_{k\in\mathbb{Z}}a_k\varphi(Nx-k)
\end{equation}
 Then the solution $\varphi$ is used in building a multiresolution for the wavelet analysis. If, for example, conditions
can be placed on (\ref{eqj1}) which yield $\ltwor$-solutions to
(\ref{eqj2}), then the closed subspace $V_0$ spanned by the
translates $\left\{\varphi(x-k)\,|\,k\in\mathbb{Z}\right\}$ is
invariant under the scaling operator
\begin{equation}\label{equ}
Uf(x)=\frac{1}{\sqrt{N}}f\left(\frac{x}{N}\right),\quad (x\in\mathbb{R})
\end{equation}
i.e. $U\left(V_0\right)\subset V_0$. Setting
$$V_j:=U^j\left(V_0\right),\quad(j\in\mathbb{Z})$$
we get the resolution
$$...V_2\subset V_1\subset V_0\subset V_{-1}\subset V_{-2}...$$
from which wavelets can be constructed as in \cite{Dau92}.
\par
The cascade operator is defined on $\ltwor$ from the masking
coefficients by:
$$M_a\psi=\sqrt{N}\sum_{n\in\mathbb{Z}}a_n\psi(N\cdot-n).$$
The scaling function $\varphi$ is then a fixed point for the
cascade operator, it satisfies the scaling equation
$$M_a\varphi=\varphi.$$
Now set
$$p(\psi_1,\psi_2)\left(e^{it}\right)=\sum_{n\in\mathbb{Z}}e^{int}\int_{\mathbb{R}}\overline{\psi_1}(x)\psi(x-n)\,dx.$$
The relation between the Ruelle operator $R_{m_0,m_0}$ and the cascade operator $M_a$ is
$$R_{m_0,m_0}\left(p(\psi_1,\psi_2)\right)=p\left(M_a\psi_1,M_a\psi_2\right),$$
and this makes the transfer operator an adequate tool in the
analysis of the orthogonality relations.
\par
One of the fundamental problems in wavelet theory is to give
necessary and sufficient conditions on $m_0$ such that the the
translates of the scaling function
$\left\{\varphi(\cdot-n)\,|\,n\in\mathbb{Z}\right\}$ form an
orthonormal set. There are two well known results that answer this
question: one due to Lawton (\cite{Law91a}), which says that one
such condition is that $R_{m_0,m_0}$, as an operator on continuous
functions, has 1 as a simple eigenvalue, the other, due to A.Cohen
(\cite{Co90}), which says that the orthogonality is equivalent to
the fact that $m_0$ has no nontrivial cycles (a cycle is a set
$\{z_1,...,z_p\}$ with $z_1^N=z_2,...,z_{p-1}^N=z_p, z_p^N=z^1$
and $|m_0(z_i)|=\sqrt{N}$ for all $i\in\{1,...,p\}$; the trivial
cycle is $\{1\}$).
\par
The peripheral spectral analysis in this paper will elucidate,
among other things, why these two conditions are equivalent.
\par
The wavelet theory gives a representation of the algebra
$\mathfrak{A}_N$ (i.e. the $C^{*}$-algebra generated by two
unitary operators $U$ and $V$ subject to the relation
$UVU^{-1}=V^N$) on $\ltwor$. $U$ is the scaling operator in
(\ref{equ}) and $V$ is the translation by 1
$$V:\psi\mapsto\psi(\cdot-1).$$
In fact we also have a representation of $\linft$ on $\ltwor$ given by
$$\pi(f)\psi=\sum_{n\in\mathbb{Z}}c_n\psi(\cdot-n),$$
for $f=\sum_{n\in\mathbb{Z}}c_nz^n\in\linft$.
\par
The scaling equation (\ref{eqj2}) can be rewritten as
$$U\varphi=\pi(m_0)\varphi.$$
This representation of $\mathfrak{A}_N$ together with the scaling
function $\varphi$ is called the wavelet representation.
\par
In \cite{Jor} it is proved that there is a 1-1 correspondence
between positive solutions to $R_{m_0,m_0}h=h$ and representations
of $\mathfrak{A}_N$. These representations are in fact given by
the unitary $U:\mathcal{H}\rightarrow\mathcal{H}$, a
representation $\pi:\linft\rightarrow\mathcal{B}(\mathcal{H})$
satisfying
$$U\pi(f)=\pi\left(f\left(z^N\right)\right)U,\quad(f\in\linft)$$
and $\varphi\in\mathcal{H}$ with
$U\varphi=\pi\left(m_0\right)\varphi$.

We reproduce here the theorem:
\begin{theorem}
\label{thjorg}\

\begin{enumerate}
\item \label{Thmax+b.3(1)}Let $m_{0}\in L^{\infty}\left(  \mathbb{T}\right)
$, and suppose $m_{0}$ does not vanish on a subset of $\mathbb{T}$ of positive
measure. Let%
\begin{equation}
\left(  Rf\right)  \left(  z\right)
=\frac{1}{N}\sum_{w^{N}=z}\left| m_{0}\left(  w\right)  \right|
^{2}f\left(  w\right) ,\qquad f\in
L^{1}\left(  \mathbb{T}\right) . \label{eqax+b.14}%
\end{equation}
Then there is a one-to-one correspondence between the data
\textup{(\ref{Thmax+b.3(1)(1)})} and \textup{(\ref{Thmax+b.3(1)(2)})} below,
where \textup{(\ref{Thmax+b.3(1)(2)})} is understood as equivalence classes
under unitary equivalence:

\begin{enumerate}
\renewcommand{\theenumi}{\relax}

\item \label{Thmax+b.3(1)(1)}$h\in L^{1}\left(  \mathbb{T}\right)  $, $h\geq
0$, and%
\begin{equation}
R\left(  h\right)  =h. \label{eqax+b.15}%
\end{equation}

\item \label{Thmax+b.3(1)(2)}$\tilde{\pi}\in\operatorname*{Rep}\left(
\mathfrak{A}_{N},\mathcal{H}\right)  $, $\varphi\in\mathcal{H}$, and the unitary
$U$ from $\tilde{\pi}$ satisfying%
\begin{equation}
U\varphi=\pi\left(  m_{0}\right)  \varphi. \label{eqax+b.16}%
\end{equation}
\end{enumerate}

\item \label{Thmax+b.3(2)}From \textup{(\ref{Thmax+b.3(1)(1)})}$\rightarrow
$\textup{(\ref{Thmax+b.3(1)(2)}),} the correspondence is given by%
\begin{equation}
\ip{\varphi}{\pi\left( f\right) \varphi}_{\mathcal{H}}=\int_{\mathbb{T}%
}fh\,d\mu, \label{eqax+b.17}%
\end{equation}
where $\mu$ denotes the normalized Haar measure on $\mathbb{T}$.

From \textup{(\ref{Thmax+b.3(1)(2)})}$\rightarrow$%
\textup{(\ref{Thmax+b.3(1)(1)}),} the correspondence is given by%
\begin{equation}
h\left(  z\right)  =h_{\varphi}\left(  z\right)  =\sum_{n\in\mathbb{Z}}%
z^{n}\ip{\pi\left( e_{n}\right) \varphi}{\varphi}_{\mathcal{H}}.
\label{eqax+b.18}%
\end{equation}

\item \label{Thmax+b.3(3)}When \textup{(\ref{Thmax+b.3(1)(1)})} is given to
hold for some $h$, and $\tilde{\pi}\in\operatorname*{Rep}\left(  \mathfrak{A}%
_{N},\mathcal{H}\right)  $ is the corresponding cyclic representation with
$U\varphi=\pi\left(  m_{0}\right)  \varphi$, then the representation is unique
from $h$ and \textup{(\ref{eqax+b.17})} up to unitary equivalence: that is, if
$\pi^{\prime}\in\operatorname*{Rep}\left(  \mathfrak{A}_{N},\mathcal{H}^{\prime
}\right)  $, $\varphi^{\prime}\in\mathcal{H}^{\prime}$ also cyclic and
satisfying%
\begin{align*}
\ip{\varphi^{\prime}}{\pi^{\prime}\left( f\right) \varphi^{\prime}} &
=\int_{\mathbb{T}}fh\,d\mu\\%
\intertext{and}%
U^{\prime}\varphi^{\prime} &  =\pi^{\prime}\left(  m_{0}\right)
\varphi^{\prime},
\end{align*}
then there is a unitary isomorphism $W$ of $\mathcal{H}$ onto $\mathcal{H}%
^{\prime}$ such that $W\pi\left(  A\right)  =\pi^{\prime}\left(  A\right)  W$,
$A\in\mathfrak{A}_{N}$, and $W\varphi=\varphi^{\prime}$.
\end{enumerate}
\end{theorem}

\begin{definition}
\label{cyclicreph} Given $h$ as in theorem \ref{thjorg} call
$\left(\pi, \mathcal{H}, \varphi\right)$ the cyclic representation
of $\mathfrak{A}_N$ associated to $h$.
\end{definition}

\par
In the case of the orthogonality of the translates of the scaling
function $\varphi$, the wavelet representation is in fact the
cyclic representation corresponding to the unique fixed point of
the Ruelle operator $R_{m_0,m_0}$, which is the constant function
$1$.
\par
We will also need the results from \cite{Dut} which show the
connection between solutions to $R_{m_0,m_0'}h=h$ and operators
that intertwine these representations. Here are those results:

\begin{theorem}
\label{thdut1} Let $m_0,m_0'\in\linft$ be non-singular and
$h,h'\in\lonet$,  $h,h'\geq0$, $R_{m_0,m_0}(h)=h$,
$R_{m_0',m_0'}(h')=h'$. Let $\left(\pi, \mathcal{H},
\varphi\right)$, $\left(\pi ', \mathcal{H}', \varphi '\right)$
be the cyclic representations corresponding to $h$ and $h'$ respectively. \\
If $h_0\in\lonet$, $R_{m_0,m_0'}\left(h_0\right)=h_0$ and $\left|
h_0\right| ^2\leq chh'$ for some $c>0$ then there exists a unique
operator $S:\mathcal{H}'\rightarrow\mathcal{H}$ such that
$$SU'=US \,,\quad S\pi '(f)=\pi (f)S\,,\quad (f\in\linft)$$
$$\ip{\varphi}{\pi (f)S\varphi '}=\int_{\mathbb{T}}fh_0\,d\mu\,,\quad(f\in\linft)$$
Moreover $\left\| S\right\|\leq\sqrt{c}$.
\end{theorem}

\begin{theorem}
\label{thdut2} Let $m_0,m_0',h,h',(\pi,\mathcal{H},\varphi),(\pi
',\mathcal{H}',\varphi')$ be as in theorem \ref{thdut1}.
\\
Suppose $S:\mathcal{H}'\rightarrow\mathcal{H}$ is a bounded operator that satisfies
$$SU'=US \,,\quad S\pi '(f)=\pi (f)S\,,\quad (f\in\linft)$$
Then there exists a unique $h_0\in\lonet$ such that
$$R_{m_0,m_0'}h_0=h_0$$
and
$$\ip{\varphi}{S\pi '(f)\varphi '}=\int_{\mathbb{T}}fh_0\,d\mu\,,\quad (f\in\linft)$$
Moreover
$$\left| h_0\right|^2\leq\left\|S\right\|^2hh'\,\mbox{ almost everywhere on }\mathbb{T}.$$
\end{theorem}
\par
We will see how each cycle of $m_0$ gives rise to a representation
of $\mathfrak{A}_N$ (hence to a positive solution for
$R_{m_0,m_0}h=h$).
\par
We will also give the concrete form for the cyclic representation
corresponding to the constant function $1$ when $m_0$ satisfies
(\ref{relm01})-(\ref{relm04}). When the wavelets are not
orthogonal (the case of tight frames), the representations become
more complicated.

\section{\label{spectral} Peripheral spectral analysis}
We begin this section by analyzing the intertwining operators a
little bit further. We will see that the commutant of the cyclic
representation associated to a positive $h$ with $R_{m_0,m_0}h=h$
is abelian and we will find the eigenfunction $h$ that corresponds
to the composition of two intertwining operators that correspond
to $h_1$ and $h_2$ respectively.
\par
 In \cite{Jor}, corollary 3.9 it is proved that the cyclic representation $(\mathcal{H}_h,\pi_h,\varphi_h)$ corresponding to some $h\geq 0$ with $R_{m_0,m_0}h=h$ is given by:
$$\mathcal{H}_h:=\left\{\left(\xi_0,...,\xi_n,...\right)\, | \, \int_{\mathbb{T}}R_{m_0,m_0}^n\left(\left|\xi_n\right|^2h\right)\, d\mu<\infty, R_{m_0,m_0}\left(\xi_{n+1}h\right)=\xi_nh \right\}$$
$$\pi_h(f)\left(\xi_0,...\xi_n,...\right)=\left(f(x)\xi_0,...,f\left(z^N\right)\xi_n,...\right)\,, \quad (f\in\linft)$$
$$U_h\left(\xi_0,...,\xi_n,...\right)=\left(m_0(z)\xi_1,...,m_0\left(z^{N^n}\right)\xi_{n+1},...\right)$$
$$\ip{\left(\xi_0,...,\xi_n,...\right)}{\left(\eta_0,...\eta_n,...\right)}=\lim_{n\rightarrow\infty}\int_{\mathbb{T}}R_{m_0,m_0}^n\left(\overline{\xi_n}\eta_nh\right)\, d\mu$$
and
$$\varphi_h=\left(1,1,...,1,...\right)$$
Also, we have the subspaces $H_0^h\subset H_1^h\subset...\subset
H_n^h\subset...\subset\mathcal{H}_h$ whose union is dense
 in $\mathcal{H}_h$ where
$$H_n^h:=\left\{\left(\xi_0,...,\xi_n,...\right)\in\mathcal{H}_h\, |\, \xi_{n+k}(z)=\xi_n\left(z^{N^k}\right)\mbox {, for } k\geq0\right\}$$
The set
$$\mathcal{V}_n^h:=\left\{U_h^{-n}\pi_h(f)\varphi_h\, |\, f\in\linft\right\}$$
is dense in $H_n^h$ for all $n\geq 0$ and $U_h^nH_n^h=H_0^h$.
\par
Some notations: if $m_0$ and $h$ are as in theorem \ref{thjorg},
then, we denote by $\left(\mathcal{H}_h,\pi_h,\varphi_h\right)$
the cyclic representation associated to $h$.
\par
If $m_0$, $m_0'$, $h, h'$ and $h_0$ are as in theorem \ref{thdut1}
then denote by $S_{h,h',h_0}$ the intertwining operator from
$\mathcal{H}_{h'}$ to $\mathcal{H}_h$ given by the aforementioned
theorem.
\par
Sometime we will omit the subscripts.
\begin{lemma}\label{lemma1}
Let $P_{H_0^h}$ be the projection onto the subspace $H_0^h$.
\\
Then $P_{H_0^h}S_{h,h',h_0}P_{H_0^{h'}}$ is multiplication by
$\frac{h_0}{h}$ on $H_0^{h'}$ i.e.
\begin{align*}P_{H_0^h}S_{h,h',h_0}P_{H_0^{h'}}\left(\xi(z),\xi\left(z^N\right),...,\xi\left(z^{N^n}\right),...\right)=\\
=\left(\xi(z)\frac{h_0(z)}{h(z)},\xi\left(z^N\right)\frac{h_0\left(z^N\right)}{h\left(z^N\right)},...,
\xi\left(z^{N^n}\right)\frac{h_0\left(z^{N^n}\right)}{h\left(z^{N^n}\right)},...\right)
\end{align*}
\end{lemma}
\begin{proof}
Denote
$S\varphi_{h'}=\left(\varphi_0^S,...,\varphi_n^S,...\right)$. Then
for all $f\in\linft$
\begin{align*}
\int_{\mathbb{T}}fh_0\, d\mu&=\ip{\left(1,1,...,1,...\right)}{\pi_h(f)\left(\varphi_0^S,...,\varphi_n^S,...\right)}\\
&=\lim_{n\rightarrow\infty}\int_{\mathbb{T}}R_{m_0,m_0}^n\left(f\left(z^{N^n}\right)\varphi_n^Sh\right)\, d\mu\\
&=\lim_{n\rightarrow\infty}\int_{\mathbb{T}}f(z)\varphi_0^Sh\,d\mu=\int_{\mathbb{T}}f\varphi_0^Sh\,d\mu
\end{align*}
thus $\varphi_0^S=\frac{h_0}{h}$. Consider again an $f\in\linft$
arbitrary.
\begin{align*}
P_{H_0^h}SP_{H_0^{h'}}\pi_{h'}(f)\varphi_{h'}&=P_{H_0^h}S\pi_{h'}(f)\varphi_{h'}=P_{H_0^h}\pi_h(f)S\varphi_{h'}=\\
&=P_{H_0^h}\left(f(z)\varphi_0^S,...,f\left(z^{N^n}\right)\varphi_n^S,...\right)=\\
&=\left(f(z)\varphi_0^S,...,f\left(z^{N^n}\right)\varphi_0^S\left(z^{N^n}\right),...\right)
\end{align*}
This calculation shows that $P_{H_0^h}SP_{H_0^{h'}}$ is
multiplication by $\frac{h_0}{h}$ on $\mathcal{V}_0^{h'}$, so, by
density, on $H_0^{h'}$.
\end{proof}
\begin{lemma}\label{lemma2}
$P_{H_n^h}S_{h,h',h_0}P_{H_n^{h'}}$ converges to $S_{h,h',h_0}$ in
the strong operator topology.
\end{lemma}
\begin{proof}
Let $\xi\in\mathcal{H}_{h'}$.
\begin{align*}
\left\| P_{H_n^h}SP_{H_n^{h'}}\xi-S\xi\right\|&\leq\left\| P_{H_n^h}SP_{H_n^{h'}}\xi-P_{H_n^h}S\xi\right\|+\left\|P_{H_n^h}S\xi-S\xi\right\|\\
&\leq
\left\|P_{H_n^h}\right\|\left\|S\right\|\left\|P_{H_n^{h'}}\xi-\xi\right\|+\left\|P_{H_n^{h}}S\xi-S\xi\right\|\rightarrow
0,\mbox{ as } n\rightarrow\infty
\end{align*}
because the subspaces $H_n^h$ form an increasing sequence whose
union in dense in $\mathcal{H}_h$ (and similarly for $H_n^{h'}$ ).
\end{proof}
\begin{theorem}\label{th1}
The commutant $\pi_h(\mathfrak{A}_n)'$ is abelian.
\end{theorem}
\begin{proof}
Consider $S_1,S_2\in\pi_h(\mathfrak{A}_n)'$. Then, according to
theorem \ref{thdut2}, $S_1=S_{h_1}$, $S_2=S_{h_2}$, for some
$h_1,h_2$ with $R_{m_0,m_0}h_i=h_i$, $\left| h_i\right|\leq c_i h$
, $i\in\{1,2\}$. Let $\xi\in\mathcal{H}_h$. It has a decomposition
$\xi=\xi_0+\eta$ with $\xi_0\in H_0^h$ and $\eta\in
{H_0^h}^{\bot}$. Using lemma \ref{lemma1}
\begin{align*}
\left(P_{H_0^h}S_1P_{H_0^h}\right)&\left(P_{H_0^h}S_2P_{H_0^h}\right)(\xi)=P_{H_0^h}S_1P_{H_0^h}S_2\xi_0=\\
&=P_{H_0^h}S_2P_{H_0^h}S_1\xi_0=\left(P_{H_0^h}S_2P_{H_0^h}\right)\left(P_{H_0^h}S_1P_{H_0^h}\right)\xi
\end{align*}
Since $P_{H_n^h}=U^{-n}P_{H_0^h}U^n$ it follows that
$P_{H_n^h}S_1P_{H_n^h}$ and $P_{H_n^h}S_2P_{H_n^h}$ also commute.
Lemma \ref{lemma2} can be used to get $S_1S_2$ as the strong limit
of
$\left(P_{H_n^h}S_1P_{H_n^h}\right)\left(P_{H_n^h}S_2P_{H_n^h}\right)$.
Similarly for $S_2S_1$. And as the limit is unique we must have
$S_1S_2=S_2S_1$.
\end{proof}
\par
Next, suppose we have two intertwining operators
$S_1:\mathcal{H}_h\rightarrow\mathcal{H}_{h'}$,
$S_2:\mathcal{H}_{h'}\rightarrow\mathcal{H}_{h''}$ which come from
$h_1$ and $h_2$ respectively. Then $S_2S_1$ is also an
intertwining operator so it must come from some $h_3$. We want to
find the relation between $h_1,h_2$ and $h_3$.
\begin{theorem}\label{th2}
If $S_{h_1}:\mathcal{H}_h\rightarrow\mathcal{H}_{h'}$ and
$S_{h_2}:\mathcal{H}_{h'}\rightarrow\mathcal{H}_{h''}$ are
intertwining operators then, if $S_{h_3}=S_{h_2}S_{h_1}$. We have
for all $f\in\linft$:
$$\int_{\mathbb{T}}|f(z)|^2R_{m_0,m_0}^n\left(\left|\frac{h_1}{h'}\frac{h_2}{h''}-\frac{h_3}{h''}\right|^2h''\right)\,d\mu\rightarrow 0$$
\end{theorem}
\begin{proof}
We begin with a calculation. For $f\in\linft$:
\begin{align*}
P_{H_n^{h'}}S_1P_{H_n^{h}}&\left(U_{h}^{-n}\pi_h(f)\varphi_h\right)=U_{h'}^{-n}P_{H_0^{h'}}U_{h'}^nS_1U_h^{-n}P_{H_0^{h}}U_n^hU_h^{-n}\pi_h(f)\varphi_h=\\
&=U_{h'}^{-n}P_{H_0^{h'}}S_1P_{H_0^h}\pi_h(f)\varphi_h\\
&=U_{h'}^{-n}\pi_{h'}\left(f\frac{h_1}{h'}\right)\varphi_{h'}
\end{align*}
For the second equality we used the fact that $S_1$ is
intertwining and for the last one, lemma \ref{lemma1}.

$$\left(P_{H_n^{h''}}S_2P_{H_n^{h'}}\right)\left(P_{H_n^{h'}}S_1P_{H_n^h}\right)\left(U_h^{-n}\pi_h(f)\varphi_h\right)=$$
$$=\left(P_{H_n^{h''}}S_2P_{H_n^{h'}}\right)U_{h'}^{-n}\pi_{h'}\left(f\frac{h_1}{h'}\right)\varphi_{h'}$$
\begin{equation}\label{rel1}
=U_{h''}^{-n}\pi_{h''}\left(f\frac{h_1}{h'}\frac{h_2}{h''}\right)\varphi_{h''}
\end{equation}
Similarly
\begin{equation}\label{rel2}
\left(P_{H_n^{h''}}S_2S_1P_{H_n^{h}}\right)\left(U_h^{-n}\pi_h(f)\varphi_h\right)=U_{h''}^{-n}\pi_h\left(f\frac{h_3}{h''}\right)\varphi_h
\end{equation}
Using (\ref{rel1}), (\ref{rel2}) and the notation
$m_0^{(n)}(z):=m_0(z)m_0\left(z^N\right)...m_0\left(z^{N^{n-1}}\right)$
, we have
$$\left\|\left(P_{H_n^{h''}}S_2P_{H_n^{h'}}\right)\left(P_{H_n^{h'}}S_1P_{H_n^h}\right)\left(\pi_h(f)\varphi_h\right)-
\left(P_{H_n^{h''}}S_2S_1P_{H_n^{h}}\right)\left(\pi_h(f)\varphi_h\right)\right\|_{\mathcal{H}_{h''}}=$$
$$\left\|\left(P_{H_n^{h''}}S_2P_{H_n^{h'}}\right)\left(P_{H_n^{h'}}S_1P_{H_n^h}\right)U_h^{-n}\pi_h\left(f\left(z^{N^n}\right)m_0^{(n)}\right)\varphi_h-\right.$$
$$\left.-\left(P_{H_n^{h''}}S_2S_1P_{H_n^{h}}\right)U_h^{-n}\pi_h\left(f\left(z^{N^n}\right)m_0^{(n)}\right)\varphi_h\right\|_{\mathcal{H}_{h''}}=$$
$$=\left\|U_{h''}^{-n}\left(\pi_{h''}\left(f\left(z^{n^n}\right)m_0^{(n)}(z)\frac{h_1}{h'}\frac{h_2}{h''}\right)\right)\varphi_{h''}-\right.$$
$$\left.-U_{h''}^{-n}\pi_{h''}\left(f\left(z^{N^n}\right)m_0^{(n)}(z)\frac{h_3}{h'}\right)\varphi_{h''}\right\|_{\mathcal{H}_{h''}}$$
$$=\int_{\mathbb{T}}\left|f\left(z^{N^n}\right)\right|^2\left|m_0^{(n)}(z)\right|^2\left|\frac{h_1}{h'}\frac{h_2}{h''}-\frac{h_3}{h''}\right|^2h''\,d\mu=$$
$$\int_{\mathbb{T}}|f(z)|^2R_{m_0,m_0}^n\left(\left|\frac{h_1}{h'}\frac{h_2}{h''}-\frac{h_3}{h''}\right|^2h''\right)\,d\mu$$
But, by lemma \ref{lemma2}, the first term in this chain of
equalities converges to 0 for all $f\in\linft$ so we obtain the
desired conclusion.
\end{proof}
\begin{corollary}\label{cor1}
If $S_{h_1},S_{h_2}\in\pi_h\left(\mathfrak{A}_N\right)'$,
$S_{h_3}=S_{h_1}S_{h_2}$ and $h\in\linft$ then
$$\int_{\mathbb{T}}|g|\left|R_{m_0,m_0}^n\left(\frac{h_1h_2}{h}\right)-h_3\right|\,d\mu\rightarrow 0,\quad (g\in\linft)$$
\end{corollary}
\begin{proof}

We will need the following inequality
\begin{equation}
\left|R_{m_0,m_0}^n(\xi h)\right|^2\leq R_{m_0,m_0}^n\left(|\xi|^2h\right)h\label{rel3}
\end{equation}
This can be proved using Schwartz's inequality:

$$\left|R_{m_0,m_0}^n(\xi h)\right|^2=\left|\frac{1}{N^n}\sum_{w^{N^n}=z}\left|m_0^{(n)}(w)\right|^2\xi(w)h(w)\right|^2\leq$$
$$\leq\left(\frac{1}{N^n}\sum_{w^{N^n}=z}\left|m_0^{(n)}(w)\right|^2|\xi(w)|^2h(w)\right)\left(\frac{1}{N^n}\sum_{w^{N^n}=z}\left|m_0^{(n)}(w)\right|^2h(w)\right)=$$
$$=R_{m_0,m_0}^n\left(|\xi|^2h\right)R_{m_0,m_0}^nh=R_{m_0,m_0}^n\left(|\xi|^2h\right)h.$$

Now take $g\in\linft$ and $f=gh^{1/2}$ in theorem \ref{th2}
$(h=h'=h'')$. We have:

$$\left(\int_{\mathbb{T}}|g|\left|R_{m_0,m_0}^n\left(\frac{h_1h_2}{h}\right)-h_3\right|\,d\mu\right)^2\leq\int_{\mathbb{T}}|g|^2\left|R_{m_0,m_0}^n\left(\frac{h_1h_2}{h}\right)-h_3\right|^2\,d\mu=$$
$$=\int_{\mathbb{T}}|g|^2\left|R_{m_0,m_0}^n\left(\left(\frac{h_1}{h}\frac{h_2}{h}-\frac{h_3}{h}\right)h\right)\right|^2\,d\mu$$
$$\leq\int_{\mathbb{T}}|g|^2hR_{m_0,m_0}^n\left(\left|\frac{h_1}{h}\frac{h_2}{h}-\frac{h_3}{h}\right|^2h\right)\,d\mu$$
$$=\int_{\mathbb{T}}|f|^2R_{m_0,m_0}^n\left(\left|\frac{h_1}{h}\frac{h_2}{h}-\frac{h_3}{h}\right|^2h\right)\,d\mu\rightarrow 0$$
\end{proof}
\par
In the sequel, we consider intertwining operators that correspond
to continuous eigenfunctions $h$. We will prove that if $h_1$ and
$h_2$ are continuous and $S_{h_3}=S_{h_1}S_{h_2}$ then $h_3$ must
be also continuous. The fundamental result needed here is from
\cite{BraJo}:
\begin{theorem}\label{th54brajo}
Let $m_0$ be a function on $\mathbb{T}$ satisfying
$m_0\in\mbox{Lip}_1(\mathbb{T})$, $R_{m_0,m_0}1\leq1$ and consider
the restriction of $R_{m_0,m_0}$ to $\mbox{Lip}_1(\mathbb{T})$
going into $\mbox{Lip}_1(\mathbb{T})$. It follows that
$R_{m_0,m_0}$ has at most a finite number
$\lambda_1,...,\lambda_p$ of eigenvalues of modulus 1,
$$\left|\lambda_i\right|=1,$$
and $R$ has a decomposition
\begin{equation}\label{eq1brajo54}
R_{m_0,m_0}=\sum_{i=1}^p\lambda_iT_{\lambda_i}+S,
\end{equation}
where $T_{\lambda_i}$ and $S$ are bounded operators from
$\mbox{Lip}_1(\mathbb{T})$ to $\mbox{Lip}_1(\mathbb{T})$,
$T_{\lambda_i}$ have  finite-dimensional range, and
\begin{equation}\label{eq2brajo54}
T_{\lambda_i}^2=T_{\lambda_i},\quad T_{\lambda_i}T_{\lambda_j}=0\mbox{ for }i\neq j,\quad T_{\lambda_i}S=ST_{\lambda_i}=0,
\end{equation}
and there exist positive constants $M$, $h$ such that
\begin{equation}
\left\|S^n\right\|_{\mbox{Lip}_1(\mathbb{T})\rightarrow\mbox{Lip}_1(\mathbb{T})}
\leq M/(1+h)^n
\end{equation}
for $n=1,2,...$. Furthermore
$\left\|R_{m_0,m_0}\right\|_{\infty\rightarrow\infty}\leq1$, and
there is a constant $M$ such that
\begin{equation}
\left\|S^n\right\|_{\infty\rightarrow\infty}\leq M
\end{equation}
for $n=1,2,...$.
\par
Finally, the operators $T_{\lambda_i}$ and $S$ extend to bounded
operators $C(\mathbb{T})\rightarrow C(\mathbb{T})$ ,  and the
properties (\ref{eq1brajo54}) and (\ref{eq2brajo54}) still hold
for this extension. Moreover
$$\lim_{n\rightarrow\infty}S^nf=0,\quad f\in C(\mathbb{T}),$$
$$T_{\lambda_i}(f)=\lim_{n\rightarrow\infty}\frac{1}{n}\sum_{k=1}^n\lambda_i^{-k}R_{m_0,m_0}^k(f),\quad f\in C(\mathbb{T}).$$
\end{theorem}
\begin{proof}
Everything is contained in \cite{BraJo}, theorem 3.4.4,
proposition 4.4.4 and its proof.
\end{proof}

\begin{theorem}\label{th3}
Assume $m_0$ is Lipschitz, $R_{m_0,m_0}1\leq1$, $h\geq 0$ is
continuous, $R_{m_0,m_0}h=h$. If
$S_{h_1},S_{h_2}\in\pi_h\left(\mathfrak{A}_N\right)'$, with
$h_1,h_2$ continuous and $S_{h_3}=S_{h_1}S_{h_2}$ then $h_3$ is
also continuous and
$$h_3=T_1\left(\frac{h_1h_2}{h}\right)=\lim_{n\rightarrow\infty}R_{m_0,m_0}^n\left(\frac{h_1h_2}{h}\right),\mbox{ uniformly.}$$
\end{theorem}
\begin{proof}
By corollary \ref{cor1} we have:
\begin{equation}\label{rel4}
\int_{\mathbb{T}}gR^n\left(\frac{h_1h_2}{h}\right)\,d\mu\rightarrow\int_{\mathbb{T}}gh_3\,d\mu\,\quad (g\in\linft)
\end{equation}
\par
Also, observe that $\frac{h_1h_2}{h}$ is continuous because
$|h_1|\leq c_1h$, $|h_2|\leq c_2h$ for some positive constants
$c_1,c_2$, and if $x_0\in\mathbb{T}$ with $h(x_0)=0$ then
$h_1(x_0)=0$, $h_2(x_0)=0$ and $\left|\frac{h_1h_2}{h}\right|\leq
c_2h_1$ Relation (\ref{rel4}) implies that for all $g\in\linft$
$$\int_{\mathbb{T}}\frac{1}{m}\sum_{n=0}^{m-1}R^n\left(\frac{h_1h_2}{h}\right)\,d\mu\rightarrow\int_{\mathbb{T}}gh_3\,d\mu$$
However, by theorem \ref{th54brajo}, we have
$$\frac{1}{m}\sum_{n=0}^{m-1} R^n\left(\frac{h_1h_2}{h}\right)\rightarrow T_1\left(\frac{h_1h_2}{h}\right),\mbox{ uniformly.}$$
Therefore $h_3=T_1\left(\frac{h_1h_2}{h}\right)$.
\par
Next we want to prove that
$R^n\left(\frac{h_1h_2}{h}\right)\rightarrow h_3$ uniformly. By
\cite{BraJo}, proposition 4.4.4, this is equivalent to
$T_{\lambda_i}\left(\frac{h_1h_2}{h}\right)=0$ for $\lambda_i\neq
1$.
\par
From (\ref{rel4}) it follows, using theorem \ref{th54brajo}, that
\begin{equation}\label{rel5}
\sum_{\lambda_1\neq 1}\lambda_i^n\int_{\mathbb{T}}gT_{\lambda_i}\left(\frac{h_1h_2}{h}\right)\,d\mu\rightarrow 0
\end{equation}
for all $g\in\linft$.
\par
But $T_{\lambda_i}\left(\frac{h_1h_2}{h}\right)$ are eigenvectors
corresponding to different eigenvalues so, some are 0 and the rest
are linearly independent. For all $i$ with
$T_{\lambda_i}\left(\frac{h_1h_2}{h}\right)\neq0$ we can find a
$g_i\in\linft$ such that
$\int_{\mathbb{T}}g_iT_{\lambda_i}\left(\frac{h_1h_2}{h}\right)\,d\mu=1$
and
$\int_{\mathbb{T}}g_iT_{\lambda_j}\left(\frac{h_1h_2}{h}\right)\,d\mu=0$
for $\lambda_j\neq\lambda_i$ (this can be obtain from the fact
that $\linft$ is the dual of $\lonet$ which contains the vectors
$T_{\lambda_i}\left(\frac{h_1h_2}{h}\right)$ ). Then, if we use
(\ref{rel5}) for $g_i$, we get that $\lambda_i^n\rightarrow 0$
whenever $T_{\lambda_i}\left(\frac{h_1h_2}{h}\right)\neq 0$,
$\lambda_i\neq 1$, which is clearly absurd unless all
$T_{\lambda_i}\left(\frac{h_1h_2}{h}\right)$ are 0, for
$\lambda_i\neq 1$. Thus, as we have mentioned before, this implies
that $R^n\left(\frac{h_1h_2}{h}\right)\rightarrow h_3$.
\end{proof}
\begin{corollary}\label{cor2}
If $h\in C(\mathbb{T})$, $h\geq0$, $R_{m_0,m_0}h=h$ then the space
$$\left\{h_0\in C(\mathbb{T})\,|\, R_{m_0,m_0}h_0=h_0, |h_0|\leq ch\right\}$$
is a finite dimensional abelian $C^*$-algebra under the pointwise
addition and multiplication by scalars, complex conjugation and
the product given by $h_1*h_2$ defined by
$S_{h_1*h_2}=S_{h_1}S_{h_2}$.
\end{corollary}
\begin{proof}
Everything follows from theorem \ref{th3} and theorem \ref{th1}.
For the finite dimensionality see \cite{BraJo} or \cite{CoRa90}.
\end{proof}
\begin{remark}\label{rem1}
When $h=1$ the $C^*$-algebra structure given in corollary
\ref{cor2} is the same as the one introduced in \cite{BraJo},
theorem 5.5.1.
\end{remark}
\par
Now we will show how each $m_0$-cycle (see definition \ref{def11}
below) gives rise to a continuous solution $h\geq 0$,
$R_{m_0,m_0}h=h$. In the end we will see that any eigenfunction
$R_{m_0,m_0}h=h$ is a linear combination of eigenfunctions coming
from such cycles.
\begin{definition}\label{def11}
Let  $m_0\in C(\mathbb{T})$. An $m_0$-cycle is a set
$\{z_1,...,z_p\}$ contained in $\mathbb{T}$ such that
$z_i^N=z_{i+1}$ for $i\in\{1,...,p-1\}$, $z_p^N=z_1$ and
$|m_0(z_i)|=\sqrt{N}$ for $i\in\{1,...,p\}$.
\end{definition}

\par
First, we consider the eigenfunction that corresponds to the cycle
$\{1\}$. This appears in many instances and it is the one that
defines the scaling function in the theory of multiresolution
approximations (see \cite{Dau92}, \cite{BraJo}).

\begin{proposition}\label{prop11}
Let $m_0\in\mbox{Lip}_1(\mathbb{T})$ with $m_0(1)=\sqrt{N}$,
$R_{m_0,m_0}1=1$. Define
$$\varphi_{m_0,1}(x)=\prod_{k=1}^{\infty}\frac{m_0\left(\frac{x}{N^k}\right)}{\sqrt{N}},\quad (x\in\mathbb{R})$$
\begin{enumerate}
\item\label{prop11i}
$\varphi_{m_0,1}$ is a well defined, continuous function and it
belongs to $\ltwor$.
\item\label{prop11ii}
If $h_{m_0,1}=\operatorname*{Per}\left|\varphi_{m_0,1}\right|^2$
is Lipschitz ( trigonometric polynomial if $m_0$ is one ), where
$$\operatorname*{Per}(f)(x):=\sum_{k\in\mathbb{Z}}f(x+2k\pi),\quad(x\in[0,2\pi],f:\mathbb{R}\rightarrow\mathbb{C}).$$

Also $R_{m_0,m_0}h_{m_0,1}=h_{m_0,1}$, $h_{m_0,1}(1)=1$,
$h_{m_0,1}$ is 0 on every $m_0$-cycle disjoint of $\{1\}$.
\item\label{prop11iii}
If $U_1:\ltwor\rightarrow\ltwor$,
$\left(U_1\xi\right)(x)=\sqrt{N}\xi(Nx)$ and
$\pi_1(f):\ltwor\rightarrow\ltwor$, $\pi_1(f)(\xi)=f\xi$ for all
$f\in\linft$, then $\left(U_1,\pi_1,\varphi_{m_0,1}\right)$ define
the cyclic representation corresponding to $h_{m_0,1}$.
\item\label{prop11iv} The commutant of the representation from (\ref{prop11iii}) is
$$\left\{M_f\,|\, f\in\linfr, f(Nx)=f(x) \mbox{ a.e. }\right\},$$
where $M_f$ is the operator of multiplication by $f$.
\item\label{prop11v}
$h_{m_0,1}$ is minimal, in the sense that if $0\leq h'\leq
ch_{m_0,1}$, $c>0$, $h'$ continuous and $R_{m_0,m_0}h'=h'$ then
there exists a $\lambda\geq0$ such that $h'=\lambda h_{m_0,1}$.
\item\label{prop11vi}
If $h\geq0$ is continuous, $R_{m_0,m_0}h=h$ and $h(1)=1$ then
$h\geq h_{m_0,1}$.
\end{enumerate}
\end{proposition}
\begin{proof}
(\ref{prop11i}) See \cite{Dau92} or \cite{BraJo}.
\par
(\ref{prop11ii}) See \cite{BraJo}, theorem 5.1.1 and lemma 5.5.6.
\par
For (\ref{prop11iv}) see \cite{Dut}. Also, in \cite{Dut} it is
proved that we are dealing with a representation of
$\mathfrak{A}_N$ (it is the Fourier transform of the wavelet
representation mentioned in the introduction). We only need to
check that $\varphi_{m_0,1}$ is cyclic for this representation.
\par
Consider $P$ the projection onto the subspace generated by
$\pi_1\left(\mathfrak{A}_N\right)\varphi_{m_0,1}$. We prove first
that $P$ commutes with the representation. Take
$A\in\pi_1\left(\mathfrak{A}_N\right)$, $A$ selfadjoint. If $B\in
\pi_1\left(\mathfrak{A}_N\right)$ then
$A(B\varphi_{m_0,1})\in\pi_1\left(\mathfrak{A}_N\right)\varphi_{m_0,1}$
so $PA(B\varphi_{m_0,1})=A(B\varphi_{m_0,1})$. So $PAP=AP$. Then
$$AP=PAP=(PAP)^*=(AP)^*=PA$$
so $P$ commutes with $A$, and since any member of
$\pi_1\left(\mathfrak{A}_N\right)$ is a linear combination of
selfadjoint operators from this set, it follows that $P$ lies in
the commutant of the representation. Then, by (\ref{prop11iv}),
$P=M_f$ for some $f\in\linfr$ with $f(Nx)=f(x)$ a.e.. As $P$ is a
projection $f^2=f=\overline f$ so $f=\chi_A$ for some subset $A$
of the real line. But $P\varphi_{m_0,1}=\varphi_{m_0,1}$ so
$\varphi_{m_0,1}\chi_A=\varphi_{m_0,1}$ a.e.. Since
$\varphi_{m_0,1}(0)=1$ and $\varphi_{m_0,1}$ is continuous, it
follows that $A$ contains a neighbourhood of $0$. This, coupled
with the fact that $\chi_A(Nx)=\chi_A(x)$ a.e., imply that
$\chi_A=1$ a.e. so $P$ is the identity and thus
$\pi_1\left(\mathfrak{A}_N\right)\varphi_{m_0,1}$ is dense, which
means exactly that $\varphi_{m_0,1}$ cyclic.

\par
(\ref{prop11v}) Consider $h'$ as mentioned in the hypothesis. Then
$h'$ induces a member of the commutant $S_{h'}$. By
(\ref{prop11iv}), $S_{h'}=M_{f_{h'}}$ for some $f_{h'}\in\linfr$
with $f_{h'}(Nx)=f_{h'}(x)$ a.e.. We have
$$\ip{\varphi_{m_0,1}}{S_{h'}\pi_1(f)\varphi_{m_0,1}}=\int_{\mathbb{T}}fh'\,d\mu,\quad (f\in\linft)$$
which implies that
$$h'=\operatorname*{Per}\left(\overline{\varphi_{m_0,1}}S_{h'}\varphi_{m_0,1}\right)=\operatorname*{Per}\left(f_{h'}\left|\varphi_{m_0,1}\right|^2\right)$$
We prove that $f_{h'}$ is continuous at 0.
\begin{equation}\label{rel6}
h'(x)=f_{h'}(x)\left|\varphi_{m_0,1}\right|^2(x)+\sum_{k\neq0}f_{h'}(x+2k\pi)\left|\varphi_{m_0,1}\right|^2(x+2k\pi)
\end{equation}
As
$$h_{m_0,1}(x)=\left|\varphi_{m_0,1}\right|^2(x)+\sum_{k\neq0}\left|\varphi_{m_0,1}\right|^2(x+2k\pi)$$
and $h_{m_0,1}(0)=\left|\varphi_{m_0,1}\right|^2(0)=1$ and
$h_{m_0,1}$, $\varphi_{m_0,1}$ are continuous, it follows that
$$\sum_{k\neq0}\left|\varphi_{m_0,1}\right|^2(x+2k\pi)\rightarrow 0\,\mbox{ as } x\rightarrow 0$$
Then, as $x\rightarrow0$,
$$\left|\sum_{k\neq0}f_{h'}(x+2k\pi)\left|\varphi_{m_0,1}\right|^2(x+2k\pi)|\right|\leq\left\|f_{h'}\right\|_{\infty}\sum_{k\neq0}\left|\varphi_{m_0,1}\right|^2(x+2k\pi)\rightarrow 0.$$
Using this in (\ref{rel6}) we obtain that
$\lim_{x\rightarrow0}f_{h'}(x)=h'(0)$. But $f_{h'}(Nx)=f_{h'}(x)$
a.e. so $f_{h'}=h'(0)$ a.e. which implies that $h'=h'(0)h_{m_0,1}$
.
\par
(\ref{prop11vi}) This is contained also in \cite{BraJo} but here
is a different proof. Consider $f\in\linft$, arbitrary. Define
$$\varphi_n(x)=f(x)\chi_{[-N^n\pi,N^n\pi]}h^{1/2}\left(\frac{x}{N^n}\right)\prod_{k=1}^n\frac{m_0\left(\frac{x}{N^k}\right)}{\sqrt{N}}$$
Clearly
$$\varphi_n(x)\rightarrow f(x)\varphi_{m_0,1},\quad (x\in\mathbb{R})$$
\begin{align*}
\int_{\mathbb{R}}\left|\varphi_n(x)\right|^2\,dx&=\int_{-N^n\pi}^{N^n\pi}|f|^2(x)h\left(\frac{x}{N^n}\right)
\prod_{k=1}^n\frac{|m_0|^2\left(\frac{x}{N^k}\right)}{\sqrt{N}}\,dx\\
&=\int_{-\pi}^{\pi}|f|^2\left(N^ny\right)h(y)\prod_{k=0}^{n-1}|m_0|^2\left(N^kx\right)\,dy\\
&=\int_{-\pi}^{\pi}R_{m_0,m_0}^n\left(h(y)|f|^2\left(N^ny\right)\right)\, dy\\
&=\int_{-\pi}^{\pi}|f|^2(y)h(y)\,dy
\end{align*}
Using Fatou's lemma we obtain:
\begin{align*}
\int_{\mathbb{R}}|f(x)|^2\left|\varphi_{m_0,1}\right|^2\,dy&=\int_{\mathbb{R}}\liminf_n\left|\varphi_n\right|^2\,dx\leq\\
&\leq\liminf_n\int_{\mathbb{R}}\left|\varphi_n(x)\right|^2\, dx\\
&=\int_{\mathbb{T}}|f|^2h\,d\mu
\end{align*}
and after periodization
$$\int_{\mathbb{T}}|f|^2h_{m_0,1}\,d\mu\leq\int_{\mathbb{T}}|f|^2h\,d\mu$$
As $f$ was arbitrary this shows that $h_{m_0,1}\leq h$.
\end{proof}

\par
Now we generalize a little bit, by considering a cycle $\{z_0\}$
where $z_0^N=z_0$.

\begin{proposition}\label{prop12}
Let $m_0\in\mbox{Lip}_1\left(\mathbb{T}\right)$,
$z_0\in\mathbb{T}$ with $z_0^N=z_0$,
$m_0\left(z_0\right)=\sqrt{N}e^{i\theta_0}$, $R_{m_0,m_0}1=1$.
Define
$$\varphi_{m_0,z_0}(x)=\prod_{k=1}^{\infty}\frac{e^{-i\theta_0}\alpha_{z_0}\left(m_0\right)\left(\frac{x}{N^k}\right)}{\sqrt{N}},\quad (x\in\mathbb{R})$$
where $\alpha_{\rho}(f)(z)=f(\rho z)$ for $z,\rho\in\mathbb{T}$
and $f\in\linft$.
\begin{enumerate}
\item\label{prop12i}
$\varphi_{m_0,z_0}$ is a well defined continuous function that
belongs to $\ltwor$.
\item\label{prop12ii}
$h_{m_0,z_0}:=\alpha_{z_0^{-1}}\left(\operatorname*{Per}\left|\varphi_{m_0,z_0}\right|^2\right)$
is Lipschitz (trigonometric polynomial if $m_0$ is one ),
$R_{m_0,m_0}h_{m_0,z_0}=h_{m_0,z_0}$, $h_{m_0,z_0}(z_0)=1$,
$h_{m_0,z_0}$ is 0 on every $m_0$-cycle disjoint of $\{z_0\}$.
\item\label{prop12iii}
If $U_{z_0}:\ltwor\rightarrow\ltwor$,
$U_{z_0}\xi=e^{i\theta_0}U_1\xi$  and
$\pi_{z_0}(f)(\xi)=\pi_1\left(\alpha_{z_0}(f)\right)(\xi)$ for
$f\in\linft$, then
$\left(U_{z_0},\pi_{z_0},\varphi_{m_0,z_0}\right)$ define the
cyclic representation corresponding to $h_{m_0,z_0}$.
\item\label{prop12iv}
The commutant of this representation is
$$\left\{M_f\,|\, f\in\linfr, f(Nx)=f(x) \mbox{ a.e. }\right\}.$$
\item\label{prop12v}
$h_{m_0,z_0}$ is minimal ( see proposition \ref{prop11} (\ref{prop11v}) ).
\item\label{prop12vi}
If $h\geq0$ is continuous , $R_{m_0,m_0}h=h$ and $h(z_0)=1$ then
$h\geq h_{m_0,z_0}$.
\end{enumerate}
\end{proposition}
\begin{proof}
Consider $m_0':=e^{-i\theta_0}\alpha_{z_0}(m_0)$. We check that
$m_0'$ satisfies the hypotheses of proposition \ref{prop11}.
Clearly $m_0'$ is Lipschitz, $m_0'(1)=\sqrt{N}$,
\begin{align*}
R_{m_0',m_0'}1 (z)&=\frac{1}{N}\sum_{w^N=z}\left|\alpha_{z_0}\left(m_0\right)\right|^2(w)=
\frac{1}{N}\sum_{w^N=z}\left|m_0\right|^2\left(z_0w\right)\\
&=\frac{1}{N}\sum_{y^N=z_0z}\left|m_0\right|^2(y)=R_{m_0,m_0}\left(z_0z\right)=1
\end{align*}
Thus we can apply proposition \ref{prop11} to $m_0'$.
\par
(\ref{prop12i}) $\varphi_{m_0,z_0}=\varphi_{m_0',1}$ and
everything follows.
\par
(\ref{prop12ii}) $h_{m_0,z_0}=\alpha_{z_0^{-1}}\left(h_{m_0',1}\right)$
\begin{align*}
R_{m_0,m_0}h_{m_0,z_0}(z)&=\frac{1}{N}\sum_{w^N=z}\left|m_0\right|^2(w)\alpha_{z_0^{-1}}\left(h_{m_0',1}\right)(w)\\
&=\frac{1}{N}\sum_{w^N=z}\left|m_0\right|^2(w)h_{m_0',1}\left(wz_0^{-1}\right)\\
&=\frac{1}{N}\sum_{y^N=zz_0^{-1}}\left|m_0\right|^2\left(yz_0\right)h_{m_0',1}(y)\\
&=R_{m_0',m_0'}h_{m_0',1}\left(zz_0^{-1}\right)=h_{m_0,z_0}(z)
\end{align*}
Also $h_{m_0,z_0}(z_0)=h_{m_0',1}(z_0^{-1}z_0)=1$ and, if $C$ is
an $m_0$-cycle disjoint of $\{z_0\}$ then $z_0^{-1}C$ is an
$m_0'$-cycle disjoint of $\{1\}$ and again proposition
\ref{prop11} applies.
\par
(\ref{prop12iii}) and (\ref{prop12iv}) can also be deduced from
proposition \ref{prop11}. The relation
$$U_{z_0}\pi_{z_0}(f)=\pi_{z_0}\left(f\left(z^N\right)\right)U_{z_0}$$
follows from the identity
$\alpha_{z_0}\left(f\left(z^N\right)\right)=\alpha_{z_0}(f)\left(z^N\right)$.
\par
(\ref{prop12v}) If $h'$ is as given, then $\alpha_{z_0}(h')$
satisfies: $0\leq\alpha_{z_0}(h')\leq
c\alpha_{z_0}\left(h_{m_0,z_0}\right)=ch_{m_0',1}$ and
$R_{m_0',m_0'}\alpha_{z_0}(h')=\alpha_{z_0}\left(R_{m_0,m_0}h'\right)=\alpha_{z_0}(h')$.
Then, by proposition \ref{prop11}, $\alpha_{z_0}(h')=\lambda
h_{m_0',1}$ for some $\lambda\geq0$ so $h'=\lambda h_{m_0,z_0}$.
\par
(\ref{prop12vi}) The argument is similar to the one used in
(\ref{prop12v}).
\end{proof}

\par
Using proposition \ref{prop12} we are now able to prove that each
$m_0$-cycle gives rise to a continuous solution for
$R_{m_0,m_0}h=h$.

\begin{proposition}\label{prop13}
Let $m_0\in\mbox{Lip}_1\left(\mathbb{T}\right)$, $R_{m_0,m_0}1=1$
and let $C=\{z_1,z_2=z_1^N,...,z_p=z_{p-1}^N\}$, $z_p^N=z_1$, be
an $m_0$-cycle, $m_0\left(z_k\right)=\sqrt{N}e^{i\theta_k}$ for
$k\in\{1,...,p\}$. Denote by $\theta_C=\theta_1+...\theta_p$.
Define
$$\varphi_{k,m_0,C}(x)=\prod_{k=1}^{\infty}\frac{e^{-i\theta_C}\alpha_{z_k}\left(m_0^{(p)}\right)\left(\frac{x}{N^{kp}}\right)}{\sqrt{N^p}},\quad (k\in\{1,...,p\})$$
\begin{enumerate}
\item\label{prop13i}
$\varphi_{k,m_0,C}$ is a well defined continuous function that
belongs to $\ltwor$.
\item\label{prop13ii}
Define
$g_{k,m_0,C}=\alpha_{z_k^{-1}}\left(\operatorname*{Per}\left|\varphi_{k,m_0,C}\right|^2\right)$
for all $k\in\{1,...,p\}$. $g_{k,m_0,C}$ is Lipschitz
(trigonometric polynomial if $m_0$ is one). Also
$$R_{m_0,m_0}^p g_{k,m_0,C}=g_{k,m_0,C},$$
$$R_{m_0,m_0}g_{k,m_0,C}=g_{k+1,m_0,C}$$
(we will use the notation mod $p$ that is $z_{p+1}=z_1$,
$g_{p+2,m_0,C}=g_{2,m_0,C}$ etc.)
\par
$g_{k,m_0,C}\left(z_j\right)=\delta_{kj}$, $g_{k,m_0,C}$ is 0 on
every $m_0$-cycle disjoint of $C$.
\item\label{prop13iii}
Define $h_{m_0,C}=\sum_{k=1}^pg_{k,m_0,C}$. Then $h_{m_0,C}$ is
Lipschitz (trigonometric polynomial if $m_0$ is one).
$R_{m_0,m_0}h_{m_0,C}=h_{m_0,C}$, $h_{m_0,C}(z_k)=1$ for all
$k\in\{1,...,p\}$ and $h_{m_0,C}$ is 0 on every $m_0$-cycle
disjoint of $C$.
\item\label{prop13iv}
$h_{m_0,C}$ is minimal.
\item\label{prop13v}
If $h\geq0$ is continuous, $R_{m_0,m_0}h=h$ and $h$ is $1$ on $C$
then $h\geq h_{m_0,C}$.
\item\label{prop13vi}
If $U_C:\ltwor^p\rightarrow\ltwor^p$,
$$U_C\left(\xi_1,...,\xi_p\right)=\left(e^{i\theta_1}U_1\xi_2,...,e^{i\theta_{p-1}}U_1\xi_p,e^{i\theta_p}U_1\xi_1\right)$$
and for $f\in\linft$, $\pi_C(f):\ltwor^p\rightarrow\ltwor^p$,
$$\pi_C(f)\left(\xi_1,...,\xi_p\right)=\left(\pi_1\left(\alpha_{z_1}(f)\right)\left(\xi_1\right),...,\pi_1\left(\alpha_{z_p}(f)\right)\left(\xi_p\right)\right)$$
then $\left(U_C,\pi_C,\left(\varphi_{1,m_0,C},...,\varphi_{p,m_0,C}\right)\right)$ is the cyclic representation \\
corresponding to $h_{m_0,C}$.
\item\label{prop13vii}
The commutant of this representation is
$$\left\{M_{f_1}\oplus...\oplus M_{f_p}\,|\, f_k\in\linfr, f_{k+1}(Nx)=f_k(x)\mbox{ a.e }, \mbox{ for } k\in\{1,...,p\}\right\}$$
\end{enumerate}
\end{proposition}

\begin{proof}
Let $m_0':=m_0^{(p)}$. Observe that
\begin{align*}
m_0'\left(z_i\right)&=m_0^{(p)}\left(z_i\right)=m_0\left(z_i\right)m_0\left(z_i^N\right)...m_0\left(z_i^{N^{p-1}}\right)=\\
&=m_0\left(z_1\right)m_0\left(z_2\right)...m_0\left(z_p\right)=\sqrt{N^p}e^{i\theta_C}
\end{align*}
\par
(\ref{prop13i}) Note that $R_{m_0^{(p)},m_0^{(p)}}=R_{m_0,m_0}^p$
so $R_{m_0^{(p)},m_0^{(p)}}1=1$. Thus (\ref{prop13i}) follows from
proposition \ref{prop12} (\ref{prop12i}) (replace $N$ by $N^p$
when working with $m_0^{(p)}$ ).
\par
(\ref{prop13ii}) If $y_1,y_2=y_1^N,...,y_q=y_{q-1}^N,y_1=y_q^N$ is
an $m_0$-cycle, then $\{y_i\}$ is an $m_0^{(p)}$-cycle. Therefore
, all assertions in (\ref{prop13ii}), except the one that relates
$g_{k,m_0,C}$ and $g_{k+1,m_0,C}$, follow from proposition
\ref{prop12} (\ref{prop12ii}).
\par
We check now (\ref{prop13vi}). $U_C$ is unitary as a composition
of unitary operators. For $f\in\linft$ we have:

$$U_C\pi_C(f)(\xi_1,...,\xi_p)=
=\left(e^{i\theta_1}\pi_1\left(\alpha_{z_2}(f)\left(z^N\right)\right)U_1\xi_2,...\right.$$
$$...\left.e^{i\theta_{p-1}}\pi_1\left(\alpha_{z_p}(f)\left(z^N\right)\right)U_1\xi_p,
e^{i\theta_p}\pi_1\left(\alpha_{z_1}(f)\left(z^N\right)\right)U_1\xi_1\right)$$
$$=\left(e^{i\theta_1}\pi_1\left(\alpha_{z_1}\left(f\left(z^N\right)\right)\right)U_1\xi_2,...
,e^{i\theta_{p-1}}\pi_1\left(\alpha_{z_{p-1}}\left(f\left(z^N\right)\right)\right)U_1\xi_p,
e^{i\theta_p}\pi_1\left(\alpha_{z_1}\left(f\left(z^N\right)\right)\right)U_1\xi_1\right)$$
$$=\pi_C\left(f\left(z^N\right)\right)U_C\left(\xi_1,...,\xi_p\right)$$
Here we used that
$\alpha_{z_{i+1}}(f)\left(z^N\right)=\alpha_{z_i}\left(f\left(z^N\right)\right)$.
\par
We must check also that
$$U_C\left(\varphi_{1,m_0,C},...,\varphi_{p,m_0,C}\right)=\pi_C\left(m_0\right)\left(\varphi_{1,m_0,C},...,\varphi_{p,m_0,C}\right)$$
To do this observe that
\begin{align*}
\alpha_{z_1}\left(m_0^{(p)}\right)(z)&=\alpha_{z_1}\left(m_0(z)\right)\alpha_{z_1}\left(m_0\left(z^N\right)\right)
\alpha_{z_1}\left(m_0\left(z^{N^{p-1}}\right)\right)=\\
&=\alpha_{z_1}\left(m_0\right)(z)\alpha_{z_2}\left(m_0\right)\left(z^N\right)...\alpha_{z_p}\left(m_0\right)\left(z^{N^{p-1}}\right)
\end{align*}
Thus
\begin{align*}
\varphi_{1,m_0,C}(x)&=\frac{e^{-i\theta_p}\alpha_{z_p}\left(m_0\right)\left(\frac{x}{N}\right)}{\sqrt{N}}
\frac{e^{-i\theta_{p-1}}\alpha_{z_{p-1}}\left(m_0\right)\left(\frac{x}{N^2}\right)}{\sqrt{N}}
...
\frac{e^{-i\theta_1}\alpha_{z_1}\left(m_0\right)\left(\frac{x}{N^p}\right)}{\sqrt{N}}\\
&\frac{e^{-i\theta_p}\alpha_{z_p}\left(m_0\right)\left(\frac{x}{N^{p+1}}\right)}{\sqrt{N}}
\frac{e^{-i\theta_{p-1}}\alpha_{z_{p-1}}\left(m_0\right)\left(\frac{x}{N^{p+2}}\right)}{\sqrt{N}}
...
\frac{e^{-i\theta_1}\alpha_{z_1}\left(m_0\right)\left(\frac{x}{N^{2p}}\right)}{\sqrt{N}}\\
...
\end{align*}
so
$$\varphi_{1,m_0,C}(x)=\prod_{k=1}^{\infty}\frac{e^{-i\theta_{1-k}}\alpha_{z_{1-k}}\left(m_0\right)\left(\frac{x}{N^k}\right)}{\sqrt{N}}$$
Similarly
$$\varphi_{i,m_0,C}(x)=\prod_{k=1}^{\infty}\frac{e^{-i\theta_{i-k}}\alpha_{z_{i-k}}\left(m_0\right)\left(\frac{x}{N^k}\right)}{\sqrt{N}}\,\quad\mbox{ for } i\in\{1,...,p\}$$
Using this formula we obtain:
\begin{align*}
U_1\varphi_{i+1,m_0,C}&=\sqrt{N}\varphi_{i+1,m_0,C}\left(Nx\right)\\
&=e^{-i\theta_i}\alpha_{z_i}\left(m_0\right)\prod_{k=2}^{\infty}\frac{e^{-i\theta_{i+1-k}}\alpha_{z_{i+1-k}}\left(m_0\right)\left(\frac{x}{N^{k-1}}\right)}{\sqrt{N}}\\
&=e^{-i\theta_i}\alpha_{z_i}\left(m_0\right)\varphi_{i,m_0,C}
\end{align*}
which shows that
$$U_C\left(\varphi_{1,m_0,C},...,\varphi_{p,m_0,C}\right)=\pi_C\left(m_0\right)\left(\varphi_{1,m_0,C},...,\varphi_{p,m_0,C}\right)$$
 Next we compute the commutant. Consider $A:\ltwor^p\rightarrow\ltwor^p$ commuting with the representation. Let $P_i$ be the
 projection onto the $i$-th component, and let $A_{ij}=P_iAP_j$. Note that
$$U_C^p\left(\xi_1,...,\xi_p\right)=\left(e^{-i\theta_C}U_1^p\xi_1,...,e^{-i\theta_C}U_1^p\xi_p\right)$$
\par
Also, since $z_i^{N^p}=z_i$, $z_i=\frac{2\pi k_i}{N^p-1}$ for some
integer $k_i$. Take any $\frac{2\pi}{N^p-1}$-periodic essentially
bounded function, $g$. Then $\alpha_{z_i}(g)=g$ so
$$\pi_C(g)\left(\xi_1,...,\xi_p\right)=\left(\pi_1(g)\xi_1,...,\pi_p(g)\xi_p\right)$$
Then $P_i$ commute with $U_C^p$ and $\pi_C(g)$ so $A_{ij}$ commute
with $U_1^p$ and $\pi_1(g)$  and, using the argument in \cite{Dut}
(proof of theorem 4.1),  (see also the proof of lemma
\ref{lemma14} below), it follows that $A_{ij}=M_{f_{ij}}$ for some
$f_{ij}\in\linfr$.
\par
Since $A$ and $\pi_C(f)$ commute for all $f\in\linft$, we have for
$i\in\{1,...,p\}$
$$\sum_{j=1}^pf_{ij}\pi_1\left(\alpha_{z_j}(f)\right)\xi_j=\pi_1\left(\alpha_{z_i}(f)\right)\sum_{j=1}^pf_{ij}\xi_j$$
Fix $k$ and take $\xi_j=0$ for all $j\neq k$, then
$$f_{ik}\pi_1\left(\alpha_{z_k}(f)\right)\xi_k=\pi_1\left(\alpha_{z_i}(f)\right)f_{ik}\xi_k$$
so $f_{ik}=0$ for $i\neq k $. Then, since $A$ commutes with $U$ we
have
$$\left(e^{i\theta_1}\sqrt{N}f_{22}(Nx)\xi_2(Nx),...,e^{i\theta_{p-1}}\sqrt{N}f_{pp}(Nx)\xi_p(Nx),e^{i\theta_p}\sqrt{N}f_{11}(Nx)\right)=$$
$$=\left(e^{i\theta_1}f_{11}(x)\sqrt{N}\xi_2(Nx),...,e^{i\theta_{p-1}}\sqrt{N}f_{p-1p-1}(x)\xi_p(Nx),e^{i\theta_p}\sqrt{N}f_{pp}(Nx)\right)$$
Therefore
$$f_{22}(Nx)=f_{11}(x)\mbox{ a.e. }$$
$$f_{33}(Nx)=f_{22}(x)\mbox{ a.e.}$$
$$\vdots$$
$$f_{11}(Nx)=f_{pp}(x)\mbox{ a.e. }$$
and (\ref{prop13vii}) follows.
\par
The cyclicity of
$\left(\varphi_{1,m_0,C},...,\varphi_{p,m_0,C}\right)$ follows as
in the proof of proposition \ref{prop11} (\ref{prop11iii}).
\par
We check that $R_{m_0,m_0}g_{i,m_0,C}=g_{i+1,m_0,C}$. Take
$f\in\linft$. We have:
\begin{align*}
\int_{\mathbb{T}}fg_{i+1,m_0,C}\,d\mu&=\ip{\varphi_{i+1,m_0,C}}{\pi_1\left(\alpha_{z_{i+1}}(f)\right)\varphi_{i+1,m_0,C}}\\
&=\ip{U_1\varphi_{i+1,m_0,C}}{U_1\pi_1\left(\alpha_{z_{i+1}}(f)\right)\varphi_{i+1,m_0,C}}=
\end{align*}
$$=\ip{e^{-i\theta_i}\pi_1\left(\alpha_{z_i}\left(m_0\right)\right)\varphi_{i,m_0,C}}{e^{-i\theta_i}\pi_1\left(\alpha_{z_{i+1}}(f)\left(z^N\right)\right)\pi_1\left(\alpha_{z_i}\left(m_0\right)\right)\varphi_{i,m_0,C}}$$
$$=\ip{\varphi_{i,m_0,C}}{\pi_1\left(\alpha_{z_i}\left(f\left(z^N\right)\right)\alpha_{z_i}\left(\left|m_0\right|^2\right)\right)\varphi_{i,m_0,C}}$$
$$=\int_{\mathbb{T}}f\left(z^N\right)\left|m_0\right|^2g_{i,m_0,C}\,d\mu=\int_{\mathbb{T}}f(z)R_{m_0,m_0}g_{i,m_0,C}\,d\mu.$$

Hence
$$R_{m_0,m_0}g_{i,m_0,C}=g_{i+1,m_0,C}.$$
\par
(\ref{prop13iii}) follows from (\ref{prop13ii}).
\par
Next we prove that $h_{m_0,C}$ is minimal. Take a continuous $h'$
with $0\leq h'\leq h_{m_0,C}$, $R_{m_0,m_0}h'=h'$. Then
$R_{m_0^{(p)},m_0^{(p)}}h'=R_{m_0,m_0}^ph'=h'$ and
$$0\leq h'\leq c\left(g_{1,m_0,C}+...+g_{p,m_0,C}\right)$$
Now we use the fact that the space
$$\left\{g\in C\left(\mathbb{T}\right)\,|\,R_{m_0^{(p)},m_0^{(p)}}g=g\right\}$$
is a $C^*$-algebra isomorphic to $C\left(\{1,...,d\}\right)$ for
some $d$ ( see corollary \ref{cor2} ), and by proposition
\ref{prop12} (\ref{prop12iv}), $g_{i,m_0,C}$ are minimal. It
follows that $h'$ can be written uniquely as
 $$h'=\alpha_1g_{1,m_0,C}+...+\alpha_pg_{p,m_0,C}$$
with $\alpha_1,...,\alpha_p\in\mathbb{C}$ ( the uniqueness comes
from the fact that $g_{i,m_0,C}$ are linearly independent, which,
in turn, is implied by (\ref{prop13ii}) ). Then
$$R_{m_0,m_0}h'=\alpha_1g_{2,m_0,c}+...+\alpha_{p-1}g_{p,m_0,C}+\alpha_pg_{1,m_0,C}$$
so, by uniqueness $\alpha_1=\alpha_2=...=\alpha_p=\alpha_1$ and
$$h'=\alpha_1\left(g_{1,m_0,C}+...+g_{p,m_0,C}\right)=\alpha_1h_{m_0,C}.$$
\par
For (\ref{prop13v}) we use a similar argument: take $h'$ as given
in the hypothesis. Then
$R_{m_0^{(p)},m_0^{(p)}}h'=R_{m_0,m_0}^ph'=h'$,
$h'\left(z_i\right)=1$ for all $i$. Using proposition \ref{prop12}
(\ref{prop12v}), we get $h'\geq g_{i,m_0,C}$ for all $i$.
\par
Now we use again the fact that $\left\{g\in
C\left(\mathbb{T}\right)\,|\,R_{m_0^{(p)},m_0^{(p)}}g=g\right\}$
is a $C^*$-algebra isomorphic to $C(\left\{1,...,d\}\right)$ and
$g_{i,m_0,C}$ are minimal, so
$$h'\geq\left(g_{1,m_0,C}+...+g_{p,m_0,C}\right)=h_{m_0,C}$$
\end{proof}

\begin{lemma}\label{lemma14}
Consider $m_0$, $m_0'$ satisfying (\ref{relm01})-(\ref{relm04}).
Let $C: z_1^N=z_2,..., z_p^N=z_1$ be an $m_0$-cycle and
$C':z_1'^N=z_2',...,z_{p'}'^N=z_1'$ be an $m_0'$-cycle,
$m_0\left(z_k\right)= \sqrt{N}e^{i\theta_k}$,
$m_0'\left(z_k'\right)=\sqrt{N}e^{i\theta_k'}$ for all $k$.
Consider the cyclic representations associated to this cycles as
in proposition \ref{prop13}, $\left(U_C,\pi_C,\varphi_C\right)$,
$\left(U_{C'},\pi_{C'},\varphi_{C'}\right)$ and let
$S:\ltwor^{p'}\rightarrow\ltwor^p$ be an intertwining operator.
Then $S=0$ if $C\neq C'$. If $C=C'$ and, after relabeling,
$z_k=z_k'$ for all $k$, $p=p'$ then, there exist
$f_1,...,f_p\in\linfr$ such that
$$S\left(\xi_1,...,\xi_p\right)=\left(f_1\xi_1,...,f_p\xi_p\right)$$
with
$$f_1(x)=e^{i\left(\theta_1-\theta_1'\right)}f_2(Nx)\mbox{, a.e.,}$$
$$\vdots$$
$$f_{p-1}(x)=e^{i\left(\theta_{p-1}-\theta_{p-1}'\right)}f_p(Nx)\mbox{, a.e.,}$$
$$f_p(x)=e^{i\left(\theta_p-\theta_p'\right)}f_1(Nx)\mbox{, a.e..}$$
\end{lemma}
\begin{proof}
Note that
$$U_C^p=e^{i\theta_C}U_1^p\oplus...\oplus e^{i\theta_C}U_1^p$$
where $\theta_C=\theta_1+...+\theta_p$. Similarly for
$U_{C'}^{p'}$. This shows that $U_C^p$ commutes with the
projections $P_i$ onto the $i$-th component.
\par
 We have $SU_{C'}^{pp'}=U_C^{pp'}S$ so $\left(P_iSP_j\right)U_{C'}^{pp'}=U_{C}^{pp'}\left(P_iSP_j\right)$,
therefore \\
$S_{ij}e^{ip\theta_C'}U_1^{pp'}=e^{ip'\theta_C}U_1^{pp'}S_{ij}$
where $S_{ij}=P_iSP_j$.
\par
Also, since $z_k^{N^p}=z_k$, $z_k$ has the form
$e^{i\frac{2l\pi}{m}}$ for all $k$ and similarly for $z_k'$. If we
take $f\in\linft$ to be $\frac{2\pi}{mm'}$-periodic, then
$\alpha_{z_k}(f)=f$, $\alpha_{z_k'}(f)=f$ for all $k$ so
$$\pi_C\left(\xi_1,...,\xi_p\right)=\left(\pi_1(f)\xi_1,...,\pi_1(f)\xi_p\right)$$
$$\pi_{C'}\left(\xi_1,...,\xi_{p'}\right)=\left(\pi_1(f)\xi_1,...,\pi_1(f)\xi_{p'}\right)$$
and again
$$S_{ij}\pi_1(f)=\pi_1(f)S_{ij}$$
Hence $S_{ij}$ commutes with $\pi_1(f)=M_f$ whenever $f\in\linfr$
is $\frac{2\pi}{mm'}$-periodic.
\par
But then also
$$\left(U_1^{-pp'}\pi_1(f)U_1^{pp'}\right)S_{ij}=S_{ij}\left(U_1^{-pp'}\pi_1(f)U_1^{pp'}\right)$$
and $U_1^{-pp'}\pi_1(f)U_1^{pp'}=M_g$ where
$g\left(N^{pp'}x\right)=f(x)$ for $x\in\mathbb{R}$ and $g$ is
$\frac{2\pi}{mm'}N^{pp'}$-periodic. By induction, it follows that
$S_{ij}$ commutes with $M_f$ whenever $f\in\linfr$ is
$\frac{2\pi}{mm'}N^{lpp'}$ periodic, $l\in\mathbb{N}$.
\par
Now take $f\in\linfr$. Define $f_l(x)=f(x)$ on
$\left[-\frac{\pi}{mm'}N^{lpp'},\frac{\pi}{mm'}N^{lpp'}\right]$
and extend it to $\mathbb{R}$ such that $f_l$ is
$\frac{2\pi}{mm'}N^{lpp'}$-periodic.
\par
We prove that $M_{f_l}$ converges to $M_f$ in the strong operator
topology. Take $\psi\in\ltwor$.
\begin{align*}
\left\|M_{f_l}\psi-M_f\psi\right\|_{\ltwor}&=\int_{\mathbb{R}}\left|f_l-f\right|^2|\psi|^2\,dx\\
&=\int_{|x|\geq\frac{\pi}{mm'}N^{lpp'}}\left|f_l-f\right|^2|\psi|^2\,dx\\
&\leq\left(2\left\|f\right\|_{\infty}^2\right)\int_{\mathbb{R}}\chi_{\left\{|x|\geq\frac{\pi}{mm'}N^{lpp'}\right\}}\left|\psi\right|^2\,dx\\
&\rightarrow 0 \mbox{, as } l\rightarrow\infty
\end{align*}
Consequently the limit holds and $M_f$ will commute also with
$S_{ij}$. As $f$ was arbitrary in $\linfr$ , using theorem IX.6.6
in \cite{Con90}, we obtain that $S_{ij}=M_{f_{ij}}$ for some
$f_{ij}\in\linfr$.
\par
Having this, we rewrite the intertwining properties. First, we
have for all $f\in\linft$:
\begin{equation}\label{eq1lemma}
\sum_{j=1}^{p'}f_{ij}\alpha_{z_j'}(f)\xi_j=\alpha_{z_i}(f)\sum_{j=1}^{p'}f_{ij}\xi_j,\quad (i\in\{1,...,p\})
\end{equation}
Fix $k\in\{1,...,p'\}$ and take $\xi_j=0$ for $j\neq k$. Then
\begin{equation}\label{eq2lemma}
f_{ik}\alpha_{z_k'}(f)\xi_k=\alpha_{z_i}(f)f_{ik}\xi_k
\end{equation}
Since $f\in\linft$ is arbitrary, it follows that $f_{ik}=0$ unless
$z_k'=z_i$.
\par
If $z_k'=z_i$ then we get $C=C'$. If $C\neq C'$ then $C\cap
C'=\emptyset$ so $f_{ij}=0$ for all $i,j$ and $S=0$.
\par
It remains to consider the case $C=C'$ and, relabeling $z_k=z_k'$
for all $k$, $p=p'$. Equation (\ref{eq2lemma}) implies that
$f_{ij}=0$ for $i\neq j$ so
$$S\left(\xi_1,...,\xi_p\right)=\left(f_1\xi_1,...,f_p\xi_p\right)$$
(we used the notation $f_i=f_{ii}$  ).
\par
The fact that $SU_{C'}=U_CS$ can be rewritten:
$$f_1(x)e^{i\theta_1'}\sqrt{N}\xi_2(Nx)=e^{i\theta_1}\sqrt{N}f_2(Nx)\xi_2(Nx)$$
$$\vdots $$
$$f_{p-1}(x)e^{i\theta_{p-1}'}\sqrt{N}\xi_p(Nx)=e^{i\theta_{p-1}}\sqrt{N}f_p(Nx)\xi_p(Nx)$$
$$f_p(x)e^{i\theta_p'}\sqrt{N}\xi_1(Nx)=e^{i\theta_p}\sqrt{N}f_1(Nx)\xi_1(Nx)$$
so
$$f_1(x)=e^{i\left(\theta_1-\theta_1'\right)}f_2(Nx)\mbox{, a.e.,}$$
$$\vdots$$
$$f_{p-1}(x)=e^{i\left(\theta_{p-1}-\theta_{p-1}'\right)}f_p(Nx)\mbox{, a.e.,}$$
$$f_p(x)=e^{i\left(\theta_p-\theta_p'\right)}f_1(Nx)\mbox{, a.e..}$$
\end{proof}

\begin{theorem}\label{th15}
Let $m_0$ satisfy (\ref{relm01})-(\ref{relm04}). Let $C_1,...,C_n$
be the $m_0$-cycles. Then, each $h\in C(\mathbb{T})$ with
$R_{m_0,m_0}h=h$ can be written uniquely as
$$h=\sum_{i=1}^n\alpha_ih_{m_0,C_i}$$
with $\alpha_i\in\mathbb{C}$. Moreover $\alpha_i=h|_{C_i}$. In
particular
$$1=\sum_{i=1}^nh_{m_0,C_i}$$
\end{theorem}
\begin{proof}
Proposition \ref{prop13} (\ref{prop13iii}) shows that
$h_{m_0,C_i}$ are linearly independent. Since the dimension of
$$\left\{h\in C(\mathbb{T}\,|\,R_{m_0,m_0}h=h\right\}$$
is $n$ ( see \cite{BraJo} ), it follows that $h_{m_0,C_i}$ form a
basis for this space. so
$$h=\sum_{i=1}^n\alpha_ih_{m_0,C_i}$$
for some $\alpha_i\in\mathbb{C}$. An application of proposition
\ref{prop13} (\ref{prop13iii}) shows that $\alpha_i=h|{C_i}$.
\end{proof}

\begin{theorem}\label{th16}
Suppose $m_0$ satisfies the conditions
(\ref{relm01})-(\ref{relm04}). Let $C_1,...,C_n$ be the
$m_0$-cycles. For each $i$ consider
$\left(U_{C_i},\pi_{C_i},\varphi_{C_i}\right)$ which give the
cyclic representation corresponding to $h_{m_0,C_i}$ ( see
proposition \ref{prop13} ). Define
$$U=U_{C_1}\oplus...\oplus U_{C_n},$$
$$\pi=\pi_{C_1}\oplus...\oplus\pi_{C_n},$$
$$\varphi=\varphi_{C_1}\oplus...\oplus\varphi_{C_n}.$$
Then $\left(U,\pi,\varphi\right)$ give the cyclic representation
corresponding to the constant function $1$. Each element $S$ in
the commutant of this representation has the form
$$S=S_{C_1}\oplus...\oplus S_{C_n}$$
where $S_{C_i}$ is in the commutant of
$\left(U_{C_i},\pi_{C_i},\varphi_{C_i}\right)$.
\end{theorem}
\begin{proof}
Since
$$1=\sum_{i=1}^nh_{m_0,C_i}$$
, for the first statement it is enough to check that $\varphi$ is
cyclic. For this we will need the commutant and then the reasoning
is the same as the one in the proofs of proposition \ref{prop11}
(\ref{prop11iii}) or proposition \ref{prop13} (\ref{prop13vi}).
But lemma \ref{lemma14} makes it clear that the elements of the
commutant have the form mentioned in the hypotesis ( see also the
proof of theorem \ref{th17} ). We also need to prove that if $S$
is in the commutant , $S=S^2=S^*$ and $S\varphi=\varphi$ then $S$
is the identity. But,
$$S=S_{C_1}\oplus...\oplus S_{C_n}$$
so $S_{C_i}=S_{C_i}^2=S_{C_i}^*$ and
$S_{C_i}\varphi_{C_i}=\varphi_{C_i}$ and, as $\varphi_{C_i}$ is
cyclic in the corresponding representation, it follows that
$S_{C_i}$ is the identity so $S=I$.
\end{proof}

\begin{theorem}\label{th17}
Suppose $m_0$ satisfies (\ref{relm01})-(\ref{relm04}). Let
$C_1,...,C_n$ be the $m_0$-cycles,
$C_i:z_{1i},z_{2i}=z_{1i}^N,...,z_{p_i\,i}=z_{p_i-1\,i}^N,z_{1i}=z_{p_i\,i}^N$
, for $i\in\{1,...,n\}$. Let $g_{k,m_0,C_i}$ be as in proposition
\ref{prop13}, $k\in\{1,...,p_i\}$, $i\in\{1,...,n\}$.
\par
If $h\in C(\mathbb{T})$, $h\neq0$ and $R_{m_0,m_0}h=\lambda h$ for
some $\lambda\in\mathbb{T}$, then there exists an
$i\in\{1,...,n\}$ such that $\lambda^{p_i}=1$ and there exist
$\alpha_i\in\mathbb{C}$, $i\in\{1,...,n\}$ such that
$$h=\sum_{i=1}^n\alpha_i\left(\sum_{k=1}^{p_i}\lambda^{-k+1}g_{k,m_0,C_i}\right)$$
and $\alpha_i=0$ if $\lambda^{p_i}\neq1$.
\end{theorem}
\begin{proof}
First note that instead of $m_0$ we can take $\left|m_0\right|$
and the problem remains the same. We have
$$\frac{1}{N}\sum_{w^N=z}\overline{\lambda m_0(w)}m_0(w)h(w)=h(z),\quad(z\in\mathbb{T})$$
so $R_{\lambda m_0,m_0}h=h$. Using theorem \ref{thdut1}, it
follows that $h$ induces an intertwining operator
$S:\mathcal{H}_{m_0}\rightarrow\mathcal{H}_{\lambda m_0}$, where
$\left(\mathcal{H}_{m_0},\pi_{m_0},\varphi_{m_0}\right)$ is the
cyclic representation corresponding to the constant function $1$
and $m_0$, and $\left(\mathcal{H}_{\lambda m_0},\pi_{\lambda
m_0},\varphi_{\lambda m_0}\right)$ is the cyclic representation
corresponding to $1$ and $\lambda m_0$.
\par
Using theorem \ref{th16} and proposition \ref{prop13}, we see that
$\mathcal{H}_{m_0}= \mathcal{H}_{\lambda m_0}$,
$\pi_{m_0}(f)=\pi_{\lambda m_0}(f)$, for $f\in\linft$,
$\varphi_{m_0}= \varphi_{\lambda m_0}$ and $U_{\lambda
m_0}=\lambda U_{m_0}$.
\par
The intertwining property of $S$ implies that
$$SU_{m_0}=\lambda U_{m_0}S$$
$$S\pi_{m_0}(f)=\pi_{m_0}(f)S,\quad (f\in\linft)$$
If $P_{C_i}$ is the projection onto the components corresponding
to the cycle $C_i$ then we see that $P_{C_i}$ commutes with both
$U_{m_0}$ and $\pi_{m_0}(f)$, for $f\in\linft$. Therefore
$$\left(P_{C_i}SP_{C_j}\right)U_{C_j}=\lambda U_{C_i}\left(P_{C_i}SP_{C_j}\right),$$
$$\left(P_{C_i}SP_{C_j}\right)\pi_{C_j}(f)=\pi_{C_i}(f)\left(P_{C_i}SP_{C_j}\right),\quad(f\in\linft)$$
Using lemma \ref{lemma14}, we obtain,
$\left(P_{C_i}SP_{C_j}\right)=0$ if $i\neq j$ and for each
$i\in\{1,...,n\}$ there exist $f_{1i},...,f_{p_i\,i}\in\linfr$
such that
$$\left(P_{C_i}SP_{C_j}\right)\left(\xi_1,...,\xi_{p_i}\right)=
\left(f_{1i}\xi_1,...,f_{p_i\,i}\xi_{p_i}\right),$$
$$f_{1i}(x)=\lambda f_{2i}(Nx)\mbox{ a.e.,}$$
$$\vdots$$
$$f_{p_i-1\,i}(x)=\lambda f_{p_i\,i}(Nx)\mbox{ a.e., }$$
$$f_{p_i\,i}(x)=\lambda f_{1i}(Nx)\mbox{ a.e.. }$$
Also, as
$$\int_{\mathbb{T}}fh\,d\mu=\ip{\varphi_{m_0}}{\pi_{m_0}(f)S\varphi_{m_0}},\quad (f\in\linft)$$
after periodization we get
$$h=\sum_{i=1}^n\sum_{k=1}^{p_i}\alpha_{z_{ki}^{-1}}\left(\operatorname*{Per}\left(f_{ki}\left|\varphi_{k,m_0,C_i}\right|^2\right)\right)$$
We want to prove that each $f_{ki}$ is continuous at $0$. Take
$i\in\{1,...,n\}$, $k\in\{1,...,p_i\}$. We know from proposition
\ref{prop13} that $g_{k,m_0,C_i}$ is 1 at $z_{ki}$ and 0 at every
other $z_{lj}$. Then
$$\left|\alpha_{z_{lj}^{-1}}\left(\operatorname*{Per}\left(f_{lj}\left|\varphi_{l,m_0,C_j}\right|^2\right)\right)\right|\leq
\left\|f_{lj}\right\|_{\infty}g_{l,m_0,C_j}$$ so this function has
limit $0$ at $z_{ki}$ for $(l,j)\neq(i,k)$. The argument used in
the proof of proposition \ref{prop11} (\ref{prop11v}) can be
repeated here to obtain that
$\lim_{x\rightarrow0}f_{ki}(x)=h(z_{ki})$.
\par
On the other hand we have
\begin{equation}\label{eq1th17}
f_{ki}\left(N^{p_i}x\right)=\lambda^{-p_i}f_{ki}(x)
\end{equation}
so if we let $x\rightarrow0$, we obtain
$h(z_{ki})=\lambda^{-p_i}h(z_{ki})$. Consequently,
$h(z_{ki})=f_{ki}=0$ or $\lambda^{p_i}=1$. Since $h\neq0$, there
exists an $i\in\{1,...,n\}$ with $\lambda^{p_i}=1$.
\par
For an $i$ with $\lambda^{p_i}\neq1$ we have $f_{ki}=0$ for all
$k\in\{1,...,p_i\}$. Now take an $i$ with $\lambda^{p_i}=1$. From
(\ref{eq1th17}) and the fact that $f_{ki}$ is continuous at $0$,
it follows that $f_{ki}$ is constant. Let $\alpha_i=f_{1i}$. Then
$f_{2i}=\lambda^{-1}\alpha_i,...,f_{p_i\,i}=\lambda^{-p_i+1}\alpha_i$
and the last assertion of the theorem is proved.
\end{proof}
\begin{corollary}\label{cor17}
Let $m_0$ as in theorem \ref{th17}. For an eigenvalue
$\lambda\in\mathbb{T}$ and $i$ with $\lambda^{p_i}=1$, define
$$h_{m_0,C_i}^{\lambda}=\sum_{k=1}^{p_i}\lambda^{-k+1}g_{k,m_0,C_i}$$. Then for each eigenvalue
$\lambda\in\mathbb{T}$, the eigenfunctions $h_{m_0,C_i}^{\lambda}$
with $\lambda^{p_i}=1$ are linearly independent. Moreover if we
define the measures
$$\nu_i^{\lambda}=\frac{1}{p_i}\sum_{k=1}^{p_i}\lambda^{k-1}\delta_{z_{ki}},\quad i\in\{1,...,n\},\lambda\in\mathbb{T},\lambda^{p_i}=1,$$
where $\delta_z$ is the Dirac measure at $z$, then
$$T_{\lambda}(f)=\sum_{i=1,\lambda^{p_i}=1}^n\nu_i^{\lambda}(f)h_{m_0,C_i}^{\lambda}.$$
\end{corollary}
\begin{proof}
First, we see that theorem \ref{th17} implies that
$h_{m_0,C_i}^{\lambda}$ with $\lambda^{p_i}=1$ span the eigenspace
 corresponding to the eigenvalue $\lambda$. Then we also note that, using proposition
\ref{prop13} (\ref{prop13ii}) we have:
\begin{equation}\label{eq1cor17}
\nu_i^{\lambda}\left(h_{m_0,C_j}^{\lambda}\right)=\delta_{ij}.
\end{equation}
This shows that $h_{m_0,C_i}^{\lambda}$ are linearly independent.
\par
On the other hand we have for all $f\in C(\mathbb{T})$, using the
fact that $C_i$ is an $m_0$-cycle:
\begin{align*}
\nu_i^{\lambda}\left(R_{m_0,m_0}(f)\right)&=\frac{1}{p_i}\sum_{k=1}^{p_i}\lambda^{k-1}\delta_{z_{ki}}\left(R_{m_0,m_0}(f)\right)=\\
&=\frac{1}{p_i}\sum_{k=1}^{p_i}\lambda^{k-1}\frac{1}{N}\sum_{w^N=z_{ki}}\left|m_0(w)\right|^2f(w)\\
&=\frac{1}{p_i}\sum_{k=1}^{p_i}\lambda^{k-1}\frac{1}{N}\left(\left|m_0\left(z_{k-1\,i}\right)\right|^2f\left(z_{k-1\,i}\right)\right.\\
&+\left.\sum_{w^N=z_{ki},w\neq z_{k-1\,i}}\left|m_0(w)\right|^2f(w)\right)\\
&=\frac{1}{p_i}\sum_{k=1}^{p_i}\lambda^{k-1}f\left(z_{k-1\,i}\right)\\
&=\lambda\nu_i^{\lambda}(f)
\end{align*}
Then, according to theorem \ref{th54brajo},
\begin{align*}
\nu_i^{\lambda}\left(T_{\lambda}(f)\right)
&=\lim_{n\rightarrow\infty}\frac{1}{n}\sum_{k=1}^n\lambda^{-k}\nu_i^{\lambda}\left(R_{m_0,m_0}^k(f)\right)\\
&=\lim_{n\rightarrow\infty}\frac{1}{n}\sum_{k=1}^n\lambda^{-k}\lambda^k\nu_i^{\lambda}(f)\\
&=\nu_i^{\lambda}(f)
\end{align*}
This, together with (\ref{eq1cor17}) and the fact that
$h_{m_0,C_i}^{\lambda}$ form a basis for the eigenspace, imply the
last equality of the corollary.
\end{proof}

\end{document}